\documentclass[11pt]{amsart}

\newcommand{\cc}{{\mathbb C}}
\newcommand{\Ff}{{\mathcal F}}
\newcommand{\pp}{{\mathbb P}}

\begin{document}

\author[C. Simpson]{Carlos Simpson}
\address{CNRS, Laboratoire J. A. Dieudonn\'e, UMR 6621
\\ Universit\'e de Nice-Sophia Antipolis\\
06108 Nice, Cedex 2, France}
\email{carlos@math.unice.fr}
\urladdr{http://math.unice.fr/$\sim$carlos/} 
\title[Formalized proof]{Formalized proof, computation, and the construction problem in algebraic geometry}


\maketitle

It has become a classical technique to turn to theoretical 
computer science to provide computational
tools for algebraic geometry. A more recent transformation is that now we also
get {\em logical} tools, and these too should be useful in the study of algebraic varieties. 
The purpose of this note is to consider a very small part of this picture, and try
to motivate the study of computer theorem-proving techniques by looking at how they might
be relevant to a particular class of problems in algebraic geometry. 
This is only an informal discussion, based more on questions and
possible research directions than on actual results.

This note amplifies
the themes discussed in my talk at the 
``Arithmetic and Differential Galois Groups'' conference 
(March 2004, Luminy), although many specific points in the discussion were only finished more
recently. 

I would like to thank: Andr\'e Hirschowitz and Marco Maggesi, for their invaluable
insights about computer-formalized mathematics as it relates to algebraic geometry and category theory;
and Benjamin Werner, M. S. Narasimhan, Alain Connes, Andy Magid and Ehud Hrushowski
for their remarks as explained below.

\section{The construction problem}

One of the basic problems we currently encounter is to give constructions of algebraic
varieties along with computations of their topological or geometric properties. We summarize here some
of the discussion in \cite{cpkg}. 

Hodge theory tells us much about what {\em cannot} happen. However, within the restrictions of Hodge
theory, we know very little about natural examples of what {\em can} happen.
While a certain array of techniques for constructing varieties is already known, these
don't yield sufficiently many examples of the complicated topological behavior we expect.
And even for the known constructions, it is very difficult to calculate the properties of
the constructed varieties. 

This has many facets. Perhaps the easiest example to state is the question of what collections
of Betti numbers (or Hodge numbers) can arise for an algebraic variety (say, smooth and maybe projective)?
For the present discussion we pass directly on to questions about the fundamental group. 
What types of $\pi _1$ can arise? We know a somewhat diverse-sounding collection
of examples: lattices, braid groups (in the quasiprojective case) \cite{MoishezonTeicher}, 
all kinds of virtually abelian groups, solvable groups \cite{SommeseVandeVen},
plenty of calculations for plane complements of line arrangements and other arrangements in low degrees 
\cite{Libgober} \cite{CohenOrlik}
\cite{AllcockCarlsonToledo}, Kodaira surfaces, many examples of non-residually
finite groups \cite{Toledo2}.
Which $\pi _1$'s have nontrivial representations?
Recall for example an old result: 
\medskip

\noindent
{\bf Theorem} \,
{\em Any nonrigid representation of a K\"ahler group in $PSL(2, \cc )$ comes by pullback from a curve.}

\medskip

Conversely, there exist nonrigid representations of rank $>2$ which don't come
by pullback from curves.  However, in a more extended sense
all of the known examples of representations come from rigid representations (which conjecturally are motivic)
and from representations on curves, by constructions involving Grothendieck's ``six operations'' (cf \cite{TMochizuki2}). 
In particular,
the irreducible components of moduli varieties of flat connexions $M_{DR}$ which are known, are all
isomorphic to moduli varieties of representations on curves. 
Nonetheless it seems likely that there are other ``new'' representations but that we don't know about them  because
it is difficult to master the computational complexity of looking for them. 

An intermediate construction might be as follows: suppose we have a family $\{ V_t\}$ of local systems on 
$X$, such that there is a closed locus $Z\subset M_{DR}(X)$ where $\mbox{dim} H^i(X_y,V_t)$ jumps for 
$t\in Z$. Then the family $\{ R^i\pi _{\ast}\} _{t\in Z}$ might be a component of the moduli space
of local systems on $Y$. Thus the whole topic of variation of differential Galois groups could lead
to some ``semi-new'' components in this way. Nonetheless, this doesn't go too far toward the 
basic question of finding cases where there are lots of representations for a general reason.

\section{Logic and calculation}

The construction problem results in a complex logical and 
computational situation, not directly amenable either to pure theoretical considerations,
or to brute-force calculation. This could open up the road to a new type of approach,
in a direction which was forseen by
the INRIA group in Rocquencourt, when they baptised their research group ``Logi-Cal''. The idea
behind this name
was that it is becoming necessary to combine logic and calculation. The origins of this requirement
lay in computer science, exemplified for example 
by the notions of ``proof-carrying code'' and verified and extracted programs. The ``Logi-Cal'' idea was very
cogently explained by Benjamin Werner in an expos\'e in Nice a few years ago, in which he 
described its possible applications to pure mathematics using the example of the four-color theorem. 
He explained that it would be good to have a proof of the four-color
theorem which combines computer verification of the theoretical details of the argument, with
the computer computations which form the heart of the proof. He said that we could
hope to have the whole thing contained in a single document verified by a single 
program.
 
It seems clear that this very nice idea should have repercussions for a much wider array of topics.
The possibility of combining logic and computation will open up new routes in algebraic geometry.
This is because there are questions such as those related to the construction 
problem above, which are susceptible neither to pure reasoning nor to pure
computation. At this conference Andy Magid mentionned an interesting case: he had
tried some time ago to compute examples of positive-dimensional representation varieties for
finitely presented groups with more relations than generators (cf \cite{AmorosBauer}
\cite{Gromov} \cite{Catanese}). He reported that the computational complexity
of the question (which depends on parameters like the number of generators, the number and length of the
relations, and the value of $n$ if we look for representations in $GL(n)$) became overwhelming even for 
very small parameters. This suggests that it would be impossible to envision a brute-force search for examples.
In the algebraic-geometric case, we might want to take concrete varieties, compute presentations for their fundamental
groups (using braid-group techniques for example) and then compute the representation spaces. Magid's remarks
suggest that a brute-force approach to this computation will not be feasible.
On the other hand, purely theoretical techniques are unlikely to answer the most interesting question in this
regard, namely: are there new or exceptional examples which are not accounted for by known theoretical reasons?
Thus the interest of looking for a mixed approach combining theory and computation. Implementation of
such an approach could be significantly enhanced by computer-formalized proof techniques 
providing an interface between theory and calculation.

Another example seen in this conference was Ehud Hrushowski's talk about algorithmic solutions to the 
problem of computing differential Galois groups. While showing that in principle there were algorithms 
to make the computation, it appeared likely that the complexity of the algorithms would be too great to
permit their direct implementation. Again, one would like to envision a  mixed approach in which theory
provides shortcuts in determination of the answers.

Of course this type of mixing has always taken place in mathematical work \cite{GrusonPeskine}.
There have also recently been advances in the use of algorithmic methods to attack problems
such as the topology of real varieties \cite{Basu} \cite{BasuPollackRoy}. 
The relevance of computerized formulation
of the theory part is that it might well permit the process to go much farther along, as it would make available
the advances in computational power to both sides of the interaction. Currently we can benefit from advanced 
computational power on the calculation side, but this can outstrip the capacity of theory to keep up. This
phenomenon was emphasized by Alain Connes in his talk (and subsequent comments) at the PQR conference in
Brussels, June 2003. He pointed out that with computer algebra programs he could come up with new identities 
which took pages and pages just to print out; and that it would be good to have tools for interpreting this new
information which often surpasses our classical human sensory capacities. It is possible that interface tools could
be of some help, but likely in the end that we would want to connect these things directly to theoretical proof
software---a step which might on some levels bypass human understanding altogether.

A related area in which it might be useful to have a mixture of theory and computation when looking for construction
results is the Hodge conjecture.  There are many concrete situations in which we expect to find certain 
algebraic cycles, but don't in general know that they exist. For example, the Lefschetz operators or Kunneth projectors
are automatically Hodge cycles.  It would be interesting to take explicit varieties and search for algebraic cycles
representing these Hodge classes. As in the search for representations, a brute-force approach would
probably run out of steam pretty fast, and it would be interesting to see what a mixed approach could attain.
A related question is the search for constructions of varieties where the Lefschetz or Kunneth operators are
topologically interesting, namely cases where the cohomology is not mostly concentrated in the middle dimension.

Finally we mention a more vague direction. In the above examples we are looking for constructions with a
certain desired topological or geometrical behavior. However, it may also be interesting to consider the 
question of what we get when we look at an arbitrary algebraic-geometric construction process or algorithm.
This type of question is related to the field of dynamical systems, and has been popularized by S. Wolfram.
There are probably many places to look for interesting processes in algebraic geometry.
Insofar as a given process produces an infinite, combinatorially arranged collection of output, it opens up questions
of asymptotic behavior, and more generally the arrangement of results with respect to measurable properties
on the output, as well as dependence on the algorithm in question. For this type of research it
would seem essential to have tools relating theoretical properties in algebraic geometry to algorithmic questions.

\section{The Bogomolov-Gieseker inequality for filtered local systems}

We go back to look more closely at the computational issues in 
constructing representations of algebraic fundamental groups. 
There are various different possible approaches: 
\newline
--construct the representations directly on a presentation of $\pi _1$;
\newline
--construct directly the connections $(E,\nabla )$ or the Higgs bundles $(E,\theta )$;
\newline
--in the quasiprojective case, construct directly parabolic bundles, logarithmic connections, or 
``filtered local systems''. 

Most work up to now on the first approach has already had the flavor of mixing computation and theory
\cite{MoishezonTeicher} \cite{PapadimaSuciu} \cite{Libgober} 
\cite{GreuelLossenShustin} \cite{DimcaNemethi} \cite{Broughton}. 
For the second and third approaches, there is a {\em Bogomolov-Gieseker inequality} lurking about.
The basic example is the classical $3c_2 - c_1^2 \geq 0$ for surfaces of general type, with
equality implying uniformization by the ball (and in particular the uniformization gives a representation
$\pi _1 \rightarrow SU(2,1)$). This was used in Livn\'e's construction \cite{Moishezon}.
Subsequent results, as is well-known,
concern stable vector bundles and extensions to the cases of Higgs and parabolic structures
\cite{Donaldson} \cite{UhlenbeckYau} \cite{Biquard} \cite{Li} \cite{LiNarasimhan} \cite{LiWang}
\cite{TMochizuki2} \cite{Nakajima} \cite{hbnc} \cite{SteerWren}.  
In all these cases, we only obtain representations in the case of equality, so it is hard
to find numerical genericity conditions which imply existence. The Bogomolov-Gieseker inequalities
come up in a fundamental way in the analysis of the quasiprojective
case, where it seems to be a problem of finding special configurations, say of divisors in the plane,
as well as special configurations of filtrations and weights to assign to the divisors, so that equality will hold
in the BGI.

We will look at one facet of the problem---the case of filtered local systems---for which 
at least the basic definitions are elementary. 
By a {\em filtration} of a vector space we shall mean a filtration indexed by
real numbers cf \cite{hbnc}. In particular $gr^{\alpha}_F$ is nonzero for only a finite number of reals $\alpha$.
A filtration can be multiplied by a positive real number $\lambda$: define $(\lambda \Ff )_{\alpha} :=
\Ff )_{\lambda ^{-1}\alpha}$.

Fix a surface $X$ with a divisor $D$ which we shall assume (at first) to have normal crossings. A 
{\em filtered local system} is a local system $L$ on $U:=X-D$ together with a filtration $\Ff_i$ at the nearby fiber
to each irreducible component of $D$. Recall that if $D_i$ is a component and $T_i$ a
tubular neighborhood of $D_i$ in $U_i$ then the nearby fiber is a fiber of the local system at a point
$P_i\in T_i$. We require that the filtration $\Ff_i$ be invariant under the monodromy 
over $T_i$. A parabolic version of the Riemann-Hilbert correspondence makes filtered local systems
correspond to parabolic logarithmic connexions (this was pointed out for curves in \cite{hbnc} and presumably
it works similarly in higher dimensions; also it was well-known in ${\mathcal D}$-module theory for the case
of integer filtrations).  We obtain the {\em Chern classes} of a filtered local system
denoted $c_i (L,\Ff )$ which could be defined as the parabolic Chern classes of the corresponding
parabolic logarithmic connexion. 
We have the following formulae. The first Chern class is given (as a cycle on $X$) by 
$$
c_1(L,\Ff ) = - \sum _ {\alpha , i} \alpha \, {\rm dim} \left( gr _{\alpha} ^{\Ff_i}(L_{P_i})
\right) \cdot D_i.
$$
The second Chern class combines a sum over intersection points $Q$ of the divisors, plus self-intersection
contributions of the components and the square of $c_1$.
For each intersection
point choose an ordering of the two associated indices and note them by $j_Q,k_Q$. Let $Q'$ denote
a point nearby to $Q$ (in the intersection of the tubular neighborhoods $T_{j_Q}$ and $T_{k_Q}$).
Define the {\em local contribution} 
$$
c_2(L,\Ff )_Q := -\sum _{\alpha , \beta}  \alpha \beta {\rm dim} \left( gr _{\alpha} ^{\Ff_{j_Q}}
gr _{\beta} ^{\Ff_{k_Q}}(L_{Q'}) \right) ,
$$
then 
$$
c_2(L,\Ff ) =\frac{1}{2}c_1(L,\Ff )^ 2 - 
\frac{1}{2} \sum _{\alpha , i} 
{\rm dim} \left( gr _{\alpha} ^{\Ff_i}(L_{P_i})
\right) \cdot \alpha ^ 2 (D_i . D_i) - 
 \sum _Q c_2(L,\Ff )_Q \cdot Q.
$$

The Chern classes allow us to define {\em stability and semistability} in the 
usual way by comparing the slope with slopes of subobjects. These conditions should be equivalent
on filtered local systems and parabolic logarithmic connexions. Finally, there should be a {\em harmonic theory}
comparing these objects with parabolic  logarithmic Higgs bundles---where T. Mochizuki's work 
\cite{TMochizuki} \cite{TMochizuki2} comes 
in.
We won't say anything about that here\footnote{A glance at his papers should convince the 
average reader of the value of having
the help of a computer to digest the argument.} except
 to say that it should lead to a 
{\em Bogomolov-Gieseker inequality} (BGI) which we describe in a conjectural way. Here I would like to thank
M. S. Narasimhan for pointing out recently that it would be good to investigate 
the BGI
for logarithmic objects. He had in mind the logarithmic Higgs bundle case, but it seems likely that all three 
cases would be interesting and the simplest to explain and think about is filtered local systems. 

The BGI would say that if $(L,\Ff )$ is a filtered local system which is semistable with
$c_1(L,\Ff )=0$ then $c_2(L,\Ff ) \geq 0$ and in case of equality we get some kind of 
pluriharmonic metric. The pluriharmonic metric should allow us to make a correspondence with parabolic Higgs
bundles and to use the transformations discussed in \cite{products} to obtain other different representations
of $\pi _1(U)$. 

The first case to look at is when $L$ is a trivial local system of rank $r$
which we denote by $\cc^r$.
It is easiest to understand the filtrations in this case, and also in this way we don't presuppose
having any representations of $\pi _1(U)$. Even in this case, if equality could be obtained in the BGI then
the transformations of \cite{products} would yield nontrivial representations. 
By tensoring with a rank one filtered local system, we can assume that the filtrations are
{\em balanced}:
$$
\sum _{\alpha} \alpha {\rm dim} (gr_{\alpha}^{\Ff}(\cc ^r)) = 0.
$$
This guarantees that the first Chern class will vanish. Now define the {\em product of two filtrations}
by
$$
\langle \Ff , {\mathcal G} \rangle :=
\sum _{\alpha , \beta} \alpha \beta {\rm dim} \left( gr _{\alpha} ^{\Ff}
gr _{\beta} ^{{\mathcal G}}(\cc ^r) \right) .
$$
In this case the second Chern class (as a number) becomes
$$
c_2(\cc ^r, \Ff ) = - \frac{1}{2}
\sum _{i,j} \langle \Ff_i , \Ff_j \rangle D_i . D_j.
$$
The stability condition is that if $V\subset \cc ^r$ is any proper subspace, then
$$
\sum _{\alpha , i} \alpha {\rm dim} (gr_{\alpha}^{\Ff_i}(V)) {\rm deg} (D_i) <0.
$$
The BGI can be stated as a theorem in this case:

\noindent
{\bf Theorem}
{\em If $\{ \Ff_i\}$ is a collection of filtrations
satisfying the stability condition, then $c_2(\cc ^r, \Ff ) \geq 0$, and if equality holds then
there are irreducible representations of $\pi _1(X-D)$.}

The theorem in this case is a consequence of what is known for parabolic vector bundles. Indeed,
the collection of filtrations also provides a parabolic structure for the trivial vector bundle
(with the same Chern classes). 
If we use small multiples $\{ \epsilon \Ff_i\}$, then the stability condition as described above
implies stability of the parabolic bundle, so the Bogomolov-Gieseker inequality (plus representations
in case of equality) for parabolic bundles \cite{Li} \cite{LiNarasimhan} 
gives the statement of the theorem.

It may be interesting to think of the minimum of $c_2(\cc ^r,\Ff )$ as some kind of measure of how
far we are from having representations of $\pi _1(U)$. We need to be more precise because scaling the filtrations
by a positive real number doesn't affect stability and it scales the second Chern class by the square. Put
$$
\| \Ff \| ^2 := \sum _ {i,\alpha} \| \alpha \| ^2 {\rm dim} (gr ^{\Ff_i} _{\alpha})\cdot deg (D_i),
$$
and 
$$
\Upsilon (X,D,r):= \min _{\{ \Ff_i \} }\frac{c_2(\cc ^r, \Ff )}{\| \Ff \| ^2}
$$
where the minimum is taken over collections of filtrations which give a nontrivial stable filtered structure
with $c_1(\cc ^r,\Ff ) = 0$ 
on the constant local system $\cc ^r$. The BGI says that $\Upsilon (X,D,r) \geq 0$ and 
in case of equality, there should\footnote{This would require proving that the minimum is attained.}
exist nontrivial representations of $\pi _1(X-D)$.

The above considerations lead to the question of how $\Upsilon (X,D,r)$ behaves for actual normal crossings
configurations on surfaces $X$. For simplicity, $(X,D)$ might come from a plane configuration after blowing up
(for example, a plane configuration with only multiple intersections, where we blow up once at each intersection point).
The first problem is computing $\Upsilon (X,D,r)$ and in particular calculating
the local contributions to the second Chern class at points which are not normal crossings
(discussed in \cite{Li}). Computation of $\Upsilon (X,D,r)$ involves searching through
the possible configurations of filtrations. Most importantly, we would like to create configurations of 
divisors $D_i$ in the plane which are interesting with respect to the invariant $\Upsilon $.

This might be algorithmic:
given some process for generating plane configurations, what are the distribution, asymptotic behavior
and other properties of the resulting numbers $\Upsilon$? 
But even before we get to infinite families of configurations, the simple problem of thoroughly 
analyzing what happens for specific configurations is a nontrivial computational problem.
Calculation of algebro-geometric and specially topological properties of plane configurations
goes back to Zariski and Hirzebruch, and much work in this direction continues
(see Teicher {\em et al} \cite{MoishezonTeicher} \cite{RobbTeicher} \ldots ). One of the main characteristics
of these computations is that they require significant amounts of reasoning.
Similarly, the computation of Donaldson
invariants has required a significant amount of theoretical work \cite{EllingsrudLePotierStromme}
\cite{OkonekTeleman}.
The problem we are proposing above, consideration of the behavior of the BGI and the
minima $\Upsilon$ in the setting of a 
configuration, will quite likely fit into the same mold.

Back to the theoretical level, it might be interesting to look at whether we could have a Gromovian
phenomenon of simply connected
varieties which look approximately non-simply-connected, which is to say that their ``isoperimetric 
inequalities'' are very bad, with relatively small loops being the boundaries only of very large homotopies.
Also whether Bogomolov-Gieseker quantities such as $\Upsilon (X,D,r)$ being small (but nonzero)
might detect it. And again, we would like to have information about the distribution of this phenomenon
in combinatorial families of varieties.

\section{The foundations of category theory}

Unfortunately, the visions sketched above contrast with the rather limited state of progress on the problem
of computer formalization of theoretical mathematics such as algebraic geometry. It is of course necessary
to give a thorough overview of the many projects working in this direction all over the globe; but this has been or is
in the process of being treated in other documents. In this note I will rather just describe the
current state of my own progress on this matter. 

There are two Coq develoments attached to the source file of the {\tt arxiv} version of the
present preprint.\footnote{Go to the {\tt arxiv} preprint's abstract page, then to ``other formats''
and download the ``Source'' format. The result is a tar archive containing the tex source file for the preprint but also 
the {\tt *.v} files in question. Compiles with {\tt v8.0} of the Coq proof assistant.}
One is a short self-contained file {\tt fmachine.v} which is a little demonstration of how  
pure computer-programming can be done entirely within the Coq environment (we don't even need Coq's 
program-extraction mechanism). The example which is treated is a forward-reasoning program for 
a miniature style of first-order logic (compare \cite{Ridge}). Programs such as this one 
itself may or may not be useful for proof-checking
in the future. The main point of interest is that we can write a program entirely within 
Coq; this might point the way for how to treat the programming side of things when we want to 
integrate computation with mathematical theory. The notion of Coq as a programming language was
mentionned by S. Karrmann on the Coq-club mailing list \cite{Karrmann}.

The other development continues with the environment described in \cite{stmc} where we 
axiomatized 
a very classical-looking ZFC within the type-theoretical environment of Coq, maintaining
access to the type-theory side of things via the realization parameter $R$.
This is based on a small set of axioms which purport to correspond to how types are implemented as sets,
following Werner's paper \cite{Werner}
---we don't give any argument other than refering to 
\cite{Werner} for why these axioms should be consistent. 

Here we build on this by adding basic category theory. Newer---slightly updated--versions 
of the files from \cite{stmc} are included with the present development
(in particular one has to use the versions included here and not the older ones).\footnote{With
this method of making public a continuing mathematical theory development project, the 
files bundled with a given preprint do not all represent new material: some are copies of previous ones
possibly with slight modifications, while others are new but even the new ones will themselves be 
recopied in the future.} We treat the notions of category, functor and natural transformation. 
We construct the category of functors between two given categories. Then we treat limits and colimits,
and give examples of categories. 

Most of what we have done here---and more---has already been done some time ago in
different contexts:
 Huet and Saibi, in Coq, in the context of ``setoids'' \cite{HuetSaibi};
several articles in Mizar \cite{BanacerekEtAl}; 
and also\footnote{There is an extensive discussion of references about mechanizing category theory
in a thread of the QED mailing list, circa 1996, in response to a question posted by
Clemens Ballarin. David Rydeheard mentions work in the systems Alf, LEGO and Coq, and work by 
Dyckhoff, Goguen, Hagino, Aczel, Cockett, Carmody and Walters,
Fleming, Gunther, Rosebrugh, Gray, Watjen and Struckmann,
Hasegawa, and Gehrke. 
Masami Hagiya mentions work of his student Takahisa Mohri.  
Ingo Dahn mentions a number of Mizar articles by 
Byl\'inski, Trybulec, Muzalewski, Bancerek, Darmochwal.
Roger B. Jones mentions some work of his own. 
Pratt mentions work by Bruckland and Walters, and 
tools for computation with finite categories by Rosebrugh.
And Amokrane Saibi mentions his work with Huet in Coq.
Evidently this list would have considerably to be expanded for work up to the present day.}
\cite{BarrWells}
\cite{CarmodyLeemingWalters} \cite{CarmodyWalters} 
\cite{FlemingGuntherRosebrugh} \cite{Gehrke} 
\cite{MatzMillerPothoffThomasValkema} \cite{Mohri} 
\cite{RydeheardBurstall} \cite{WatjenStruckmann}.
We don't actually claim that our present treatment has any
particular advantages over the other ones; the reason for doing it is that we hope 
it will furnish
a solid foundation for future attempts to treat a wider range of mathematical theories. 

We use the following approach to defining the notion of category. A category
is an uplet (with entries named over strings using the file \verb}notation.v} as was explained in
\cite{stmc}) consisting of the set of objects, the set of morphisms, the composition function,
the identity function, and a fifth place called the ``structure'' which is a hook allowing us
to add in additional structure in the future if called for (e.g. monoidal categories will
have the tensor product operation encoded here; closed model categories will have the 
fibration, cofibration and equivalence sets encoded here; sites will have the Grothendieck topology
encoded here etc.). The elements of the set of morphisms are themselves assumed to be
``arrows'' which are triplets having a ``source'', a ``target'' and an ``arrow'' (to take care of the 
information about the morphism). In particular, the functions \verb}source} and \verb}target} 
don't depend on which category we are in. 
Functors and natural transformations are themselves arrows, so the functions \verb}source}
and \verb}target} do a lot of work. 

We treat limits in detail, and colimits by dualizing limits. The main technical work is 
directed toward the formalized proof of the following standard theorem. 

\medskip

\noindent
{\bf Theorem}\,
{\em If \verb}a}, \verb}b}, and \verb}c} 
are categories such that \verb}b} admits limits over \verb}c}, then \verb}functor_cat a b} also admits limits 
over \verb}c}.}

\medskip

The proof is done in the file \verb}fc_limits.v}. Intricacy comes from the need to 
use  the universal property of the pointwise limits in
order to construct the structural morphisms for the limiting functor, and then
further work is needed to  show that the functor constructed in this way is actually a limit.
The corresponding result for colimits is
obtained almost immediately by dualizing---the only subtlety being that \verb}opp (functor_cat a b)}
is not equal but only isomorphic to \verb}functor_cat (opp a) (opp b)}. Because of this we need to 
make a preliminary study of the invariance of limiting properties under isomorphisms of categories.
This discussion will have to be amplified in the future when we are able to treat equivalences of categories.

The importance of this theorem
is its corollary that presheaf categories admit limits and
colimits. This will (in the future) be essential to theories of sheaves and hence topoi; and theories of
closed model categories, because many useful closed model categories take presheaf categories as their
underlying categories, and one of the main conditions for a closed model category is that it 
should admit (at least finite) limits and colimits. 

One task which is worth mentionning is that we construct examples of categories by various different 
methods, in the file \verb}cat_examples.v}. The methods include subcategories of other categories;
defining a category by its object set together with the set of arrows between each pair of objects; and
function categories which come in two flavors, depending on whether we look at functions between the
objects as sets themselves or functions between their underlying sets (denoted \verb}U x}).

A different approach is called for when we want to construct and manipulate finite categories---important
for example in relating classical limit constructions such as equalizers and fiber products, to the
notion of limit as defined in general (done for (co)equalizers in \verb}equalizer.v}
and (co)fiber-products in \verb}fiprod.v}). It doesn't seem efficient to manipulate finite sets
by directly constructing them, but instead to build them with Coq's inductive type construction and then bring them
into play using the realization parameter. 
This allows us to list the elements of a 
finite type by name, and then to manipulate them with the \verb}match} construction.
To bridge from here to the notion of category, we need to discuss the construction of categories (also
functors and natural transformations) starting from type-theoretic data: these constructions \verb}catyd},
\verb}funtyd} and \verb}nttyd} occupy a large place in  \verb}little_cat.v}.

We finish by pointing out how a theoretical category-theory development such as presented here,
is relevant for some of the more long-range projects discussed in the beginning. This discussion 
is very related to L. Chicli's thesis \cite{Chicli} in which he used Huet-Saibi's category theory
as the basis for the definition and construction of affine schemes. The basic point is that to manipulate
the fundamental objects of modern (algebraic, analytic or even differential) geometry, we need to know
what a ringed space is, and better yet a ringed site or ringed topos. 
Thus we need a theory of sheaves, and in particular a well-developed category theory,
with functor categories, limits and colimits, etc.  The next items which need to be treated in the present
development are
equivalences of categories, adjoint functors (and even fancier things like Kan extensions), over-categories,
monomorphisms and epimorphisms, then sheaves and topoi.

If we want to access more recent
developments in geometry, it will be essential to have good theories of (possibly monoidal)
closed model categories
starting with the small-object argument. On a somewhat different plane, it is clear that to
manipulate many of the geometric questions discussed above, we will need to have a good 
development of linear algebra. This presents a number of categoric aspects, for example in the
notions of additive and abelian categories (again possibly with tensor structures). 

There remain some thorny notational dilemmas still to be worked out before we can do all of this.
One example is that the right notion of ``presheaf'' is probably slightly different from that of 
a functor: we probably don't want to include the data of the target category. This is because
in general the target category will be a big category for the universe we want to work in,
whereas we would also like our presheaves to be elements of the universe, and indeed 
we don't necessarily want to specify which universe it is for a given presheaf. So we will probably
have to define a presheaf as being a modified version of a functor where the \verb}target} element
of the arrow triple is set by default to \verb}emptyset}. This is the kind of thing which is easy
to say in a few phrases, but which in practice requires writing a whole new file containing 
material similar (but not identical) to what is in \verb}functor.v}.

It seems likely 
that once the definitional work is finished, subsequent geometrical manipulations
of these objects should be fairly easy to take care of, compared with the amount of foundational work 
necessary just to give the definitions. Unfortunately, as best as I know nobody has
gotten far enough to test this out.  

\section{Finite categories}

The work on formalization of category theory, {\em a priori} a
waypoint along the path to formalizing algebraic geometry, also 
suggests its own research directions. When we are forced to look very closely
at the foundational details of a subject, there stand out certain questions which
would otherwise be overlooked in the usual rush to get on with the abstract theory.
An example, 
strongly representative of the general problem of relating theory and computation,
is the classification of
finite categories \cite{Tilson}.  For a given
finite integer $N$, how many categories are there
with $N$ morphisms? What do diagrams or other standard categorical 
constructions (functors, natural transformations,
limits, adjoints, Kan extensions \ldots ) look like in these categories, 
perhaps in terms of asymptotic behavior
with respect to $N$ but also maybe just for small fixed values? What additional structures can
these categories have? 

The question of classification of finite categories has been treated in 
\cite{Tilson} \cite{Steinberg} \cite{SteinbergTilson}
\cite{Jones} \cite{Kientzle} 
from a universal-algebra point of view.
Their idea is to define notions of {\em variety} or {\em pseudovariety} which are collections
of objects closed under
direct product and subquotient,\footnote{
These notions might be modifiable so as to be relevant to the problem of classifying representations
of algebraic fundamental groups.}
and from these references we know a lot about the structure and
classification of pseudovarieties of finite 
categories.
For example, Tilson proves a classification theorem for
locally trivial categories, those being the ones with only identity
endomorphisms of each object: the answer is that they are subquotients
of products of the two-arrow category whose limits are equalizers \cite{Tilson}.
Related are
\cite{AdamekBorceuxLackRosicky}, \cite{AlmeidaWeil} \cite{AuingerSteinberg}
\cite{JonesPustejovsky} \cite{Pin} \cite{PinPinguetWeil} \cite{Stallings} \cite{Rhodes}.
And \cite{RowellStongWang}, \cite{EtingofOstrik} discuss a similar question of classification of
finite tensor categories (but the word ``finite'' has a slightly different meaning there).

One might also ask more detailed
questions about finite categories
which are not invariant under the process of taking subquotients, and we get
a situation in many ways analogous to the algebraic-geometric questions
discussed above, leading among other things to the question of how to construct finite categories
having given properties. We can also think of further questions by analogy with the algebraic-geometric
ones. For example, the analogue of the the moduli space $M_B$ could be defined
as follows. If $\Gamma$ is a finite (or even finitely presented)
category, define the moduli stack ${\mathcal M}_B(\Gamma )$ as
the stack associated to the prestack of functors
$$
{\mathcal M}_B(\Gamma )^{\rm pre}(A) := \underline{\mbox{Hom}}(\Gamma , {\bf Mod}^{\rm proj}_A)
$$
where ${\bf Mod}^{\rm proj}_A$ is the category of projective $A$-modules.  This could have variants where we look at 
all $A$-modules or even $U$-coherent sheaves on $Spec (A)$ in the
sense of \cite{HirschowitzUcoh}. There would also be $n$-stack versions where we look at maps
into stacks of complexes or other things (and indeed we could fix any $\infty$-stack ${\mathcal G}$ 
and look at $\underline{\mbox{Hom}}(\Gamma , {\mathcal G})$). It isn't our purpose to get into the details
of this type of construction here but just to note that these should exist. We can hope in some cases to 
get geometric stacks---for example the $1$-stack ${\mathcal M}_B(\Gamma )$ as defined above is Artin-algebraic (or more
precisely its $1$-groupoid interior is algebraic). 
We can also hope that these stacks have natural open substacks with
coarse moduli varieties which could be denoted 
generically by $M_B(\Gamma )$.

Invariants of these moduli
varieties (to start with, their dimensions and irreducible components \ldots) 
would become invariants of the finite category, and we would like to know something
about their distribution, bounds, etc., and also whether we can construct finite categories
such that the moduli varieties have given behavior.  In the case when $\Gamma$ is a finitely presented
category which is free over a graph,
$M_B(\Gamma )$ is the same thing as the moduli space of {\em quivers}, and
in general the moduli space
will be a subspace of the space of quivers on the arrows of the category,
so there is already a big theory
about this (and we can expect semistability for quivers to lead to the
open substack required above). It is certainly also related to work by Lusztig, MacPherson and 
Vilonen and others on combinatorial descriptions of
perverse sheaves \cite{Lusztig}
\cite{MacPhersonVilonen} \cite{GelfandMacPhersonVilonen} \cite{Vilonen}
\cite{BradenGrinberg}. 
Which finite categories arise as specialization categories for
stratifications (and particularly naturally arising stratifications)? 
We can also ask which varieties arise as moduli spaces $M_B(\Gamma )$: this might be
relevant as a process for constructing algebraic varieties. 
 
These and any number of similar questions of differing levels of difficulty 
might provide a good proving ground for
tools combining theory and calculation.


\begin{thebibliography}{MM}


\bibitem{AdamekBorceuxLackRosicky}
J. Ad\'amek, F. Borceux, S. Lack, and J. Rosick\'y.
A classification of accessible categories. {\sc Max Kelly volume},
{\em J. Pure Appl. Alg.} {\bf 75} (2002), 7-30.

\bibitem{AdamekJohnstoneMakowskyRosicky}
J. Ad\'amek, P. Johnstone, J. Makowsky, J. Rosick\'y. 
Finitary sketches. {\em J. Symbolic Logic} {\bf 62}
(1997), 699-707.
 
\bibitem{AdamekRosicky}
J. Ad\'amek, J. Rosick\'y. On geometric and finitary sketches. 
{\sc European colloquium of category theory (Tours 1994)},
{\em Appl. Cat. Struct.} {\bf 4}
(1996), 227-240.

\bibitem{AllcockCarlsonToledo}
D. Allcock, J. Carlson, D. Toledo. Orthogonal complex hyperplane arrangements. 
{\em Symposium in Honor of C. H. Clemens (Salt Lake City 2000)}, {\sc Contemp. Math.}
{\bf 312} (2002), 1-8. 


\bibitem{AlmeidaWeil}
J. Almeida and P. Weil. Profinite categories and semidirect products.
{\em J. Pure Appl. Alg.} {\bf 123} (1998), 1-50.



\bibitem{AmorosBauer}
J. Amor\'os, I. Bauer. On the number of defining relations for nonfibered K\"ahler groups. 
{\em Internat. J. Math.} {\bf 11} (2000), 285-290.

\bibitem{AmorosBurgerEtAl}
J. Amor\'os, M. Burger, K. Corlette, D. Kotschick, D. Toledo. {\em Fundamental groups of compact K\"ahler manifolds}.
{\sc A.M.S. Mathematical Surveys and Monographs}, {\bf 44} (1996).



\bibitem{ArapuraBresslerRamachandran}
D. Arapura, P. Bressler, M. Ramachandran.
On the fundamental group of a compact K\"ahler manifold. 
{\em Duke Math. J.} {\bf 68} (1992), 477-488.

\bibitem{AuingerSteinberg}
K. Auinger, B. Steinberg.
The geometry of profinite graphs with applications to free groups and finite monoids 
{\em Trans. Amer. Math. Soc.} {\bf 356} (2004), 805-851. 


\bibitem{AurouxKatzarkov}
D. Auroux, L. Katzarkov. Branched coverings of $\cc \pp ^2$ and invariants of symplectic $4$-manifolds.
{\em Inventiones} {\bf 142} (2000), 631-673.


\bibitem{BalajiBiswasNagaraj}
V. Balaji, I. Biswas, D.  Nagaraj.
Principal bundles over projective manifolds with parabolic structure over a divisor. 
{\em Tohoku Math. J.} {\bf 53} (2001), 337-367.

\bibitem{BanacerekEtAl}
G. Bancerek, C. Byl\'inski, A. Trybulec, {\em et al}, Numerous articles on category theory in MIZAR,
{\em Journal of Formalized Mathematics}


\bibitem{BarrWells}
M. Barr, C. Wells, {\em Category theory for computer science}, Prentice Hall, New York, 1990. 

\bibitem{Basu}
S. Basu.
Different bounds on the different Betti numbers of semi-algebraic sets. 
{\sc ACM Symposium on Computational Geometry (Medford, MA, 2001)}.
{\em Discrete Comput. Geom.} {\bf 30} (2003), 65-85.

\bibitem{BasuPollackRoy}
S. Basu, R. Pollack, M.-F. Roy.
{\em Algorithms in real algebraic geometry.}
{\sc Algorithms and Computation in Mathematics} {\bf 10}
Springer-Verlag (2003). 



\bibitem{Beke}
T. Beke. Isoperimetric inequalities and the Friedlander-Milnor conjecture,
  to appear in Crelle's Journal. 

\bibitem{Bihan}
F. Bihan. Asymptotiques de nombres de Betti d'hypersurfaces projectives r\'eelles. Preprint {\tt math.AG/0312259}.


\bibitem{Biquard}
O. Biquard. Sur les fibr\'es paraboliques sur une surface complexe.
{\em J. London Math. Soc.} {\bf  53} (1996), 302-316.

\bibitem{BjornerLovaszYao}
A. Bj\"orner, L. Lov\'asz, A. C. C. Yao. 
Linear decision trees: volume estimates and topological bounds.
{\em Proceedings of the twenty-fourth annual ACM symposium on Theory of computing (Victoria, British Columbia)}
ACM Press (1992), 170-177.

\bibitem{Boalch}
P. Boalch. The Klein solution to Painleve's sixth equation.
Preprint {\tt math.AG/0308221}.

\bibitem{BodenYokogawa}
H. Boden, K. Yokogawa. Moduli spaces of parabolic Higgs bundles and parabolic $K(D)$ pairs over smooth curves, I.
{\em Internat. J. Math.} {\bf 7} (1996), 573-598.

\bibitem{BogomolovKatzarkov}
F. Bogomolov, L. Katzarkov. 
Complex projective surfaces and infinite groups. {\em Geom. Funct. Anal.} {\bf 8} (1998), 243-272.

\bibitem{BradenGrinberg}
T. Braden, M. Grinberg.
Perverse sheaves on rank stratifications.
{\em Duke Math. J.} {\bf 96} (1999), 317-362.


\bibitem{Broughton}
S. Broughton. On the topology of polynomial hypersurfaces. 
{\em Singularities, Part 1 (Arcata, Calif., 1981)} 
{\sc AMS Proc. Sympos. Pure Math.} {\bf 40} (1983),  167-178, 

\bibitem{CarmodyLeemingWalters}
S. Carmody, M. Leeming.
Walters, R. F. C. The Todd-Coxeter procedure and left Kan extensions.  {\em 
J. Symbolic Comput.} {\bf  19}  (1995),  no. 5, 459-488.
 
\bibitem{CarmodyWalters} 
S. Carmody, R. Walters. Computing quotients of actions of a free category.  
{\em Category theory (Como, 1990)}, {\sc Lecture Notes in Math.} {\bf 1488} (1991),  63-78. 


\bibitem{Catanese}
F. Catanese.  Fundamental groups with few relations. {\em Higher Dimensional Complex Varieties, Trento 1994}. 
De Gruyter (1996), 163-165.

\bibitem{Chicli}
L. Chicli. Sur la formalisation des math\'ematiques dans le Calcul des Constructions Inductives.
Thesis, Universit\'e de Nice-Sophia Antipolis (Nov. 2003). 
\newline
\verb}http://www-sop.inria.fr/lemme/Laurent.Chicli/these_chicli.ps}


\bibitem{CohenOrlik}
D. Cohen, P. Orlik. Arrangements and local systems. {\em Math. Res. Lett.} {\bf 7} (2000), 299-316.


\bibitem{DelgadoMargolisSteinberg}
M. Delgado, S.  Margolis, B. Steinberg.
Combinatorial group theory, inverse monoids, automata, and global semigroup theory. 
{\em Geometric and Combinatorial Methods in Group Theory and Semigroup Theory (Lincoln, Nebraska, 2000)}.
{\sc Internat. J. Algebra Comput.} {\bf 12} (2002), 179-211.


\bibitem{DimcaNemethi}
A. Dimca, A. Nemethi. On the monodromy of complex polynomials. {\em Duke Math. J.} {\bf 108} (2001), 199-209.

\bibitem{Donaldson}
S. Donaldson.
Anti self-dual Yang-Mills connections over complex algebraic surfaces and stable vector bundles.
{\em Proc. London Math. Soc.} {\bf 50} (1985),  1-26. 


\bibitem{Eilenberg}
S. Eilenberg. {\em Automata, Languages and Machines}. 
Academic Press, New York, Vol. A (1974), Vol B (1976).
 

\bibitem{EllingsrudLePotierStromme}
G. Ellingsrud, J. Le Potier, S. A. Stromme. Some Donaldson invariants of $\cc \pp ^2$. 
{\em Moduli of vector bundles (Sanda, 
Kyoto, 1994)} {\sc Lect. Notes Pure Appl. Math.} {\bf 179} (1996), 33-38.


\bibitem{EsnaultNoriSrinivas}
H. Esnault, M. Nori, V. Srinivas. 
Hodge type of projective varieties of low degree. 
{\em Math. Ann.} {\bf 293} (1992), 1-6.

\bibitem{EtingofOstrik}
P. Etingof, V. Ostrik.
Finite tensor categories.  Preprint {\tt math/0301027}.

\bibitem{FlemingGuntherRosebrugh}
M. Fleming, R. Gunther, R. Rosebrugh.
A database of categories. 
{\em J. Symbolic Comput.} {\bf 35} (2003), 127-135.



\bibitem{FriedmanMorgan}
R. Friedman,  J. Morgan.  Algebraic surfaces and Seiberg-Witten invariants. 
{\em J. Algebraic Geom.} {\bf 6} (1997), 445-479.

\bibitem{Gehrke}
W. Gehrke. Rewriting Techniques applied to Basic Category Theory. Technical Report 94-56, 
Research Institute for Symbolic Computation, Linz, Austria, 1994. 

\bibitem{GehrkePfalzgraf}
W. Gehrke, J.  Pfalzgraf.
Computer-aided construction of finite geometric spaces: 
automated verification of geometric constraints.
{\em J. Automat. Reason.} {\bf 26} (2001), 139-160.

\bibitem{GelfandMacPhersonVilonen}
S. Gelfand, R. MacPherson, K.  Vilonen.
Perverse sheaves and quivers. 
{\em Duke Math. J.} {\bf 83} (1996), 621-643.

\bibitem{GreuelLossenShustin}
G.-M. Greuel, C. Lossen, E.  Shustin. Plane curves of minimal degree with prescribed singularities. 
{\em Inventiones} {\bf 133} (1998), 539-580.

\bibitem{Gromov}
M. Gromov. 
Sur le groupe fondamental d'une vari\'et\'e k\"ahl\'erienne. 
{\em C. R. Acad. Sci. Paris S\'er. I Math.} {\bf 308} (1989), 67-70.

\bibitem{Gromov2}
M. Gromov.
Isoperimetry of waists and concentration of maps.
{\em Geom. Funct. Anal.} {\bf 13} (2003),  178-215.

\bibitem{GrusonPeskine}
L. Gruson, C. Peskine. Genre des courbes de l'espace projectif. 
{\em Algebraic geometry (Tromso, 1977)}, Springer {\sc L.N.M.} {\bf 687} (1978), 31-59.

\bibitem{Hasegawa}
Masahito Hasegawa.
Decomposing typed lambda calculus into a couple of categorical programming languages.
{\em Category theory and computer science (Cambridge 1995)}. 
Springer {\sc Lecture Notes in Comput. Sci.} {\bf 953} (1995), 200-219.

\bibitem{HirschowitzUcoh}
A. Hirschowitz. 
Coh\'erence et dualit\'e sur le gros site de Zariski. 
{\em Algebraic curves and projective geometry (Trento, 1988)}
{\sc Lecture Notes in Math.} {\bf 1389} (1989), 91-102.

\bibitem{HuetSaibi}
G. Huet, A. Saibi. 
Constructive category theory.
{\em Proof, language, and interaction}
{\sc Found. Comput. Ser.}
MIT Press, Cambridge (2000),  239-275.

\bibitem{InabaIwasakiSaito}
M. Inaba, K. Iwasaki, Masa-Hiko Saito.
Moduli of Stable Parabolic Connections, Riemann-Hilbert correspondence and Geometry of Painlev\'e 
equation of type VI, Part I.
Preprint {\tt math.AG/0309342}.

\bibitem{JohnstoneWraith}
P. Johnstone,  G.  Wraith, {\it Algebraic theories in toposes}, 
Lecture Notes in Mathematics, vol. 661, Springer-Verlag, 1978, pp. 141--242. 


\bibitem{Jones}
P. R. Jones,  Profinite categories, 
implicit operations and pseudovarieties of categories, J. Pure Appl. Algebra  109 (1996), 61--95.


\bibitem{JonesPustejovsky}
P. R. Jones, S. Pustejovsky. A kernel for relational morphisms of categories. 
{\sc Semigroups and Applications},
J. Howie, W. Munn, H. Weinert eds., World Scientific (1992), 152-161.

\bibitem{KapovichMillson}
M. Kapovich, J. Millson.
On representation varieties of Artin groups, projective arrangements and the fundamental 
groups of smooth projective varieties. {\em Publ. Math. I.H.E.S.} {\bf 88} (1998), 5-95.

\bibitem{Karrmann}
S. Karrmann. IO and Coq. Message on the {\tt coq-club} mailing list, Aug. 30, 2004. 

\bibitem{KatisSabadiniWalters}
P. Katis, N.  Sabadini, R. Walters. 
Bicategories of processes. 
{\em J. Pure Appl. Algebra} {\bf 115} (1997),  141-178.

\bibitem{Kientzle}
T. Kientzle. Categorical generalizations of classical monoid theory. 
Doctoral dissertation, Berkeley (1992).


\bibitem{Konno}
H. Konno. Construction of the moduli space of stable parabolic Higgs bundles on a Riemann surface. {\em J. Math.
Soc. Japan} {\bf 45} (1993), 253-276.

\bibitem{KrohnRhodes}
K. Krohn, J. Rhodes. Algebraic theory of machines, I. 
{\em Trans. Amer. Math. Soc.} {\bf 116} (1965), 450-464.

\bibitem{KulkarniWood}
R. Kulkarni, J. Wood. 
Topology of nonsingular complex hypersurfaces. {\em Adv. in Math.} {\bf 35} (1980), 239-263.

\bibitem{Li}
Jiayu Li. Hermitian-Einstein metrics and Chern number inequalities on parabolic stable bundles over K\"ahler 
manifolds.  {\em Comm. Anal. Geom.} {\bf 8} (2000), 445-475.



\bibitem{LiNarasimhan}
J. Li, M.S. Narasimhan. Hermitian-Einstein metrics on parabolic stable bundles. 
{\em Acta Math. Sin.} {\bf 15} (1999), 93-114.

\bibitem{LiNarasimhan2}
Jiayu Li, M. S. Narasimhan. A note on Hermitian-Einstein metrics on parabolic stable bundles. 
{\em Acta Math. Sinica} {\bf 17} (2001), 77-80.


\bibitem{LiWang}
J. Li, Y.-D. Wang. Existence of Hermitian-Einstein metrics on stable Higgs bundles over open K\"ahler manifolds.  
{\em Internat. J. Math.} {\bf 10} (1999), 1037-1052.

\bibitem{Libgober}
A. Libgober. Alexander polynomial of plane algebraic curves and cyclic multiple planes. 
{\em Duke Math. J.} {\bf 49} (1982), 833-851.


\bibitem{Libgober2}
A. Libgober. First order deformations for rank one local systems with a non-vanishing cohomology. 
{\em Arrangements in Boston: a Conference on Hyperplane Arrangements (1999).} 
{\sc Topology Appl.} {\bf 118} (2002), 159-168. 

\bibitem{Lusztig}
G. Lusztig.
Quivers, perverse sheaves, and quantized enveloping algebras.
{\em J. Amer. Math. Soc.} {\bf 4} (1991), 365-421.

\bibitem{MacPhersonVilonen}
R. MacPherson, K. Vilonen.
Elementary construction of perverse sheaves.
{\em Invent. Math.} {\bf 84} (1986), 403-435.

\bibitem{MakkaiReyes}
M. Makkai, G. Reyes, {\em First order categorical logic},
{\sc Lecture Notes in Mathematics} {\bf 611} (1977). 


\bibitem{MatzMillerPothoffThomasValkema}
O. Matz, A. Miller, A. Pothoff, W. Thomas, and E. Valkema.
Report on the program AMoRE. {\sc Tech. Rep.} {\bf 9507}, 
Christian Albrechts Universitat, Kiel (1994). 


\bibitem{TMochizuki}
T. Mochizuki. Asymptotic behaviour of tame nilpotent harmonic bundles with trivial parabolic structure. 
{\em J. Diff. Geom.} {\bf 62} (2002), 351-559.


\bibitem{TMochizuki2}
T. Mochizuki.
Asymptotic behaviour of tame harmonic bundles and an application to pure twistor $D$-modules.
Preprint {\tt math.DG/0312230}. 

\bibitem{Mohri}
Takahisa Mohri. 
On formalization of bicategory theory. 
{\em Theorem proving in higher order logics (Murray Hill, NJ, 1997)},
Springer
{\sc Lecture Notes in Comput. Sci.} {\bf 1275} (1997), 199-214. 

\bibitem{Moishezon}
B. Moishezon.
{\em Complex surfaces and connected sums of complex projective planes}. 
(appendix by R. Livne). Springer {\sc L.N.M.} {\bf 603} (1977). 

\bibitem{MoishezonTeicher}
B. Moishezon, M. Teicher. Braid group technique in complex geometry. I. Line arrangements in $\cc \pp ^2$.
{\em Braids (Santa Cruz, 1986)}, {\sc Contemp. Math.} {\bf 78}, A.M.S. (1988), 425-555.



\bibitem{Nakajima}
H. Nakajima. Hyper-K\"ahler structures on moduli spaces of parabolic Higgs bundles on Riemann surfaces. 
{\em Moduli of vector bundles (Sanda, 
Kyoto, 1994)} {\sc Lect. Notes Pure Appl. Math.} {\bf 179} (1996), 199-208.



\bibitem{NiRen}
L. Ni, H. Ren. Hermitian-Einstein metrics for vector bundles on complete K\"ahler manifolds. 
{\em Trans. Amer. Math. Soc.} {\bf 353} (2001), 441-456.

\bibitem{Nitsure}
N. Nitsure. Moduli of semistable logarithmic connections. {\em J. Amer. Math. Soc.} {\bf 6} (1993), 597-609. 

\bibitem{Nori}
M. Nori. Zariski's conjecture and related problems. {\em Ann. Sci. E.N.S.} {\bf 16} (1983), 305-344. 


\bibitem{OkonekTeleman}
C. Okonek, A. Teleman. Gauge theoretical equivariant Gromov-Witten invariants and the full Seiberg-Witten
invariants of ruled surfaces. {\em Comm. Math. Phys.} {\bf 227} (2002), 551-585.

\bibitem{PapadimaSuciu}
S. Papadima, A. Suciu. Higher homotopy groups of complements of complex hyperplane arrangements. 
{\em Adv. Math.} {\bf 165} (2002), 71-100.


\bibitem{Pin}
J.-E. Pin. BG = PG, a success story, 
{\em Semigroups, Formal Languages and Groups}, J. B. Fountain, ed., Kluwer (1995), 33-47. 

\bibitem{PinPinguetWeil}
J. Pin, A. Pinguet, P. Weil. Ordered categories and ordered semigroups. {\em Comm. Algebra} {\bf 30} (2002),
5651-5675.

\bibitem{Reznikov}
A. Reznikov. 
The structure of K\"ahler groups. I. Second cohomology. {\em Motives, polylogarithms and 
Hodge theory, Irvine, 1998}, 
{\sc Int. Press Lect. Ser.} {\bf 3} (II) (2002), 717-730.

\bibitem{Rhodes}
J. Rhodes. Undecidability, automata, and pseudovarieties of finite semigroups, 
{\em Int. J. Alg. Comput.} {\bf 9} (1999), 455-473.

\bibitem{Ridge}
T. Ridge. 
A Mechanically Verified, Efficient, Sound and Complete Theorem Prover For First Order Logic.
{\sc Isabelle} {\em Archive of Formal Proofs} (2004-09-28). 

\bibitem{RobbTeicher} 
A. Robb, M. Teicher.
Applications of braid group techniques to the decomposition of moduli spaces, new examples. 
{\em Special issue on braid groups and related topics (Jerusalem, 1995)}. 
{\sc Topology Appl.} {\bf 78} (1997), 143-151.

\bibitem{RosebrughSabadiniWalters}
R. Rosebrugh, N.  Sabadini, R. Walters.
Minimal realization in bicategories of automata. 
{\em Math. Structures Comput. Sci.} {\bf 8} (1998),  93-116.


\bibitem{RowellStongWang}
E. Rowell, R. Stong, Z. Wang.
Towards a classification of modular tensor categories. Preprint
\verb}www.math.ksu.edu/~zlin/jehconf/rowell.pdf}

\bibitem{RydeheardBurstall}
D. Rydeheard and R. Burstall. {\em 
Computational Category Theory}. {\sc International Series in Computer Science.} Prentice Hall, 1988.

\bibitem{hbnc}
C. Simpson. Harmonic bundles on noncompact curves. {\em J. Amer. Math. Soc.} {\bf 3}  (1990),  713-770.

\bibitem{products}
C. Simpson.
Products of matrices.
{\em Differential geometry, global analysis, and topology (Halifax, 1990)} 
{\sc CMS Conf. Proc.} {\bf 12},
Amer. Math. Soc. (1991), 157-185.

\bibitem{cpkg} 
C. Simpson. The construction problem in K\"ahler geometry. 
{\em Different Faces of Geometry}, 
M.  Gromov, S. Donaldson, Y. Eliashberg, eds.,
{\sc International Mathematical Series} {\bf 3} (2004).

\bibitem{stmc}
C. Simpson. Set-theoretical mathematics in Coq. Preprint {\tt math.LO/0402336}. 
   


\bibitem{SommeseVandeVen}
A. Sommese, A. Van de Ven. Homotopy groups of pullbacks of varieties. {\em Nagoya Math. J.} {\bf 102} (1986), 79-90.

\bibitem{SommeseVerschelde} 
A. Sommese, J. Verschelde. 
Numerical homotopies to compute generic points on positive dimensional algebraic sets. 
{\em Complexity theory, real machines, and homotopy (Oxford, 1999)}.  
{\sc J. Complexity} {\bf 16} (2000), 572-602.


\bibitem{Stallings}
J. Stallings. Topology of finite graphs. {\em Invent. Math.} {\bf 71} (1983), 551-565. 

  
\bibitem{SteerWren}
B. Steer, A. Wren. The Donaldson-Hitchin-Kobayashi correspondence for 
parabolic bundles over orbifold surfaces.
{\em Canad. J. Math.} {\bf 53} (2001), 1309-1339.


\bibitem{Steinberg}
B. Steinberg.
Semidirect products of categories and applications. 
{\em J. Pure Appl. Algebra} {\bf 142} (1999), 153-182.


\bibitem{Steinberg2}
B. Steinberg.
Fundamental groups, inverse Schützenberger automata, and monoid presentations. 
{\em Comm. Algebra} {\bf 28} (2000), 5235-5253.

\bibitem{Steinberg3}
B. Steinberg.
Finite state automata: a geometric approach.
{\em Trans. Amer. Math. Soc.} {\bf 353} (2001), 3409-3464. 

\bibitem{SteinbergTilson}
B. Steinberg, B.  Tilson. 
Categories as algebra, II. 
{\em Internat. J. Algebra Comput.} {\bf 13} (2003), 627-703.

\bibitem{Thaddeus}
M. Thaddeus. Variation of moduli of parabolic Higgs bundles. {\em J. Reine Angew. Math.} {\bf 547} (2002), 1-14. 

\bibitem{Tilson}
B. Tilson,  Categories as algebra: An essential ingredient in 
the theory of monoids, {\em J. Pure Appl. Alg.} {\bf  48} (1987), 83-198.


\bibitem{Toledo}
D. Toledo. Examples of fundamental groups of compact K\"ahler manifolds. 
{\em Bull. London Math. Soc.} {\bf 22} (1990), 339-343.


\bibitem{Toledo2}
D. Toledo.  Projective varieties with non-residually finite fundamental group. 
{\em Publ. Math. I.H.E.S.} {\bf 77} (1993), 103-119.



\bibitem{UhlenbeckYau}
K. Uhlenbeck, S.-T. Yau. 
On the existence of Hermitian-Yang-Mills connections in stable vector bundles. 
Frontiers of the mathematical sciences: 1985 (New York, 1985). 
Comm. Pure Appl. Math. 39 (1986), no. S, suppl., S257--S293.    

 
\bibitem{Vilonen}
K. Vilonen.
Perverse sheaves and finite-dimensional algebras.
{\em Trans. Amer. Math. Soc.} {\bf 341} (1994), 665-676.


\bibitem{WatjenStruckmann}
D. W\"atjen, W Struckmann.
An algorithm for verifying equations of morphisms in a category.
{\em Inform. Process. Lett.} {\bf 14} (1982), 104-108.

\bibitem{Wells}
C. Wells. Wreath product decompositions of categories, I and II. {\em Acta Sci. Math.}{\bf 52} (1988), 307-324.


\bibitem{Wells2}
C. Wells. A Krohn-Rhodes theorem for categories. {\em J. of Algebra} {\bf 64} (1980),  37-45.

\bibitem{Werner}
B. Werner. Sets in types, types in sets. 
{\em Theoretical aspects of computer software (Sendai, 1997)} {\sc Lecture Notes in Comput. Sci.} {\bf 1281} (1997),  530-546,


\bibitem{Yanofsky}
N. Yanofsky.
Coherence, homotopy and 2-theories. 
{\em $K$-Theory} {\bf 23} (2001), 203-235.

\bibitem{Yao}
Andrew C.-C. Yao.  Decision tree complexity and Betti numbers.
{\em Proceedings of the twenty-sixth annual ACM 
symposium on Theory of computing (Montreal)}  (1994), 615-624.

\bibitem{Yokogawa}
K. Yokogawa. Compactification of moduli of parabolic sheaves and moduli of parabolic Higgs sheaves.
{\em J. Math. Kyoto Univ.} {\bf 33} (1993), 451-504.


\end{thebibliography}
\end{document}